\documentclass[12pt]{article}
\usepackage{ustars}

\usepackage{amsthm,amsmath,amssymb,mathtools,multicol}

\usepackage{graphicx}

\usepackage[colorlinks=true,citecolor=black,linkcolor=black,urlcolor=blue]{hyperref}

		\usepackage{float}

\theoremstyle{plain}
\newtheorem{theorem}{Theorem}
\newtheorem{lemma}[theorem]{Lemma}

\theoremstyle{definition}
\newtheorem{definition}[theorem]{Definition}
\newtheorem{example}[theorem]{Example}

\newtheorem{question}[theorem]{Question}

\theoremstyle{remark}
\newtheorem{remark}[theorem]{Remark}

\title{\bf Distinct monomial orders with same induced orderings}

\author{Gabriel Sosa \\
\small Department of Mathematics \\[-0.8ex]
\small Purdue University \\[-0.8ex] 
\small West Lafayette, IN, USA\\
\small\tt gsosa@purdue.edu\\
}

\date{\dateline{July 1, 2014}{August 1,2014}{}}

\begin{document}

\maketitle


\begin{abstract}
We prove that the lexicographic, degree lexicographic and the degree reverse lexicographic orders for monomials in $R_n=K[X_1,\dots, X_n]$ are uniquely determined by their induced orderings, (i.e. their restrictions to $R_{n,i}=K[X_1,\dots, \hat{X_i}, \dots, X_n]$), when $n \geq 4$. We also show that for any $n \geq 4$ there are monomial orders that are not uniquely determined by their induced orderings, and provide examples of these orders for each $n$.
\end{abstract}

\section{Introduction}
Monomial orderings play a central role in computational commutative algebra, computational algebraic geometry and combinatorial commutative algebra because of the theory of Gr\"oebner basis and their applications. The properties of the classic monomial orderings (lexicographic, graded lexicographic and reverse lexicographic) have allowed for many insights on Hilbert functions, Betti numbers and regularity of monomial ideals in statements regarding Lex-segments and Generic Initial Ideals. The books \cite{CLO1}, \cite{CLO2} and \cite{KR} provide great introductions to these topics. Other results and advanced applications can be found in \cite{HH} and \cite{MS}. 

The graded lexicographic order is uniquely determined by its induced orderings when $n>3$, a fact that is equivalent to a characterization of compressed ideals due to Jeff Mermin, \cite[Theorem 3.12]{Me}. The property of being uniquely determined by induced orderings when $n>3$ is also shared by the lexicographic and reverse lexicographic orders, we present proofs in Theorem \ref{lexprop}. These facts prompted A. Conca \cite{Conca} to ask whether there is $k \in \mathbb{N}$, such that any monomial order on $R_n$, $n>k$, is uniquely determined by its induced orderings.  We answer Conca's question, Question \ref{question}, negatively in Theorem \ref{main}, and present examples of different monomial orderings with same induced orderings explicitly.

\section{Monomial orders and their matrix representations}

Let $K$ be a field, $S_n=\{X_1,\dots,X_n \}$ and $R_n=K[S_n]$ be the polynomial ring in $n$ variables, a monomial in $R_n$ is of the form $X^{\underline{\alpha}}=X_1^{\alpha_1}X_2^{\alpha_2}\cdots X_n^{\alpha_n}$. We will abuse notation and depending on the context denote the monomial  $X^{\underline{\alpha}}$ by either $\underline{\alpha}=(\alpha_1, \alpha_2, \dots, \alpha_n)$ or $\underline{\alpha}=\begin{bmatrix} \alpha_1 & \alpha_2 & \dots & \alpha_n \end{bmatrix}^T$. The set of all monomials belonging to a set $S$ will be denoted by $\text{Mon}(S)$. 

Additionally we will denote the set of $m \times n$ matrices with entries in the set $S$ by $M_{m \times n}(S)$ and the set of nonnegative integers, $\{0,1,2, \dots \}$, by $\mathbb{N}^*$.

\definition \label{monomial}{A \textbf{monomial order} in $R_n$ is a total order $\mathbf{<_{\tau}}$ in the elements of $\text{Mon}(R_n)$ with the following two additional properties:

\begin{enumerate}

\item If $u \neq 1$ then $1 <_{\tau} u$ for all $u \in \text{Mon}(R_n)$ (or equivalently $(0,\dots, 0) <_{\tau} \underline{\alpha}$ for all $\underline{\alpha} \in (\mathbb{N}^*)^n$ with $\underline{\alpha} \neq (0, \dots, 0)$ ).

\item If $u,v \in \text{Mon}(R_n)$ are such that $u <_{\tau} v$ then $u+w <_{\tau} v+w$ for all $w \in \text{Mon}(R_n)$ (or equivalently if  $\underline{\alpha} <_{\tau} \underline{\beta}$ then $\underline{\alpha}+\underline{\gamma} <_{\tau} \underline{\beta}+\underline{\gamma}$).
\end{enumerate}}

We will now give an overview of the three classical examples of monomial orders, for an in-depth study refer to \cite[Chapter 2]{CLO1}, \cite[Chapter 1]{CLO2}  or \cite[Chapter 1]{KR}.

\example{The \textbf{lexicographic order}, denoted $\boldsymbol{<_{lex}}$. We say that $\underline{\alpha}<_{lex} \underline{\beta}$ if the leftmost nonzero entry of $\underline{\beta}-\underline{\alpha}$ is positive. 

For example: $(2,2,2)<_{lex}(2,3,0)<_{lex}(3,0,3)$.

\

The \textbf{graded lexicographic order}, denoted $\boldsymbol{<_{deglex}}$. We say that $\underline{\alpha}<_{deglex} \underline{\beta}$ if $\displaystyle \sum_{i=1}^n \alpha_i<\sum_{i=1}^n \beta_i$, or, if $\displaystyle \sum_{i=1}^n \alpha_i=\sum_{i=1}^n \beta_i$ and $\underline{\alpha}<_{lex} \underline{\beta}$. 

For example: $(2,3,0)<_{deglex}(2,2,2)<_{deglex}(3,0,3)$.

\

The \textbf{reverse lexicographic order}, denoted $\boldsymbol{<_{revlex}}$. We say that $\underline{\alpha}<_{deglex} \underline{\beta}$ if $\displaystyle \sum_{i=1}^n \alpha_i<\sum_{i=1}^n \beta_i$, or, if $\displaystyle \sum_{i=1}^n \alpha_i=\sum_{i=1}^n \beta_i$ and the rightmost nonzero entry of $\underline{\beta}-\underline{\alpha}$ is negative. 

For example: $(2,3,0)<_{revlex}(3,0,3)<_{revlex}(2,2,2)$.}

\

The graded lexicography and reverse lexicographic orders are examples of \textbf{graded monomial orders} (i.e,  $\underline{\alpha}<_{\tau} \underline{\beta}$ whenever $\displaystyle \sum_{i=1}^n \alpha_i<\sum_{i=1}^n \beta_i$). 
 
{\remark \label{spread} Due to the definition of monomial order it is easy to see that $R_1$ admits only one monomial order, namely $1< x< x^2< \dots$. So in this case the lexicographic, graded lexicographic and reverse lexicographic orders coincide in $R_1$.

In $R_2$ there are infinitely many monomial orders, but $R_2$ admits only two graded monomial orders, as a consequence of Lemma \ref{lower}. In this case the graded lexicographic and reverse lexicographic orders actually coincide in $R_2$, while the lexicographic order does not coincide with them anymore. 

In $R_3$ all three classical orders are different, as it can be seen by the different ways in which the elements of the set $\{(2,3,0), (2,2,2), (3,0,3) \}$ were ordered before. \par}

\

These observations seem to indicate that as $n$ increases the classical monomial orders are growing further apart. To formalize what we mean by growing further apart we will use the concept of induced ordering.

\definition{ Let $S_{n,i}=S_n-\{x_i \}$ and $R_{n,i}=K[S_{n,i}]$. The $\boldsymbol{i^{th}}$ \textbf{induced ordering of} $\boldsymbol{<_{\tau}}$ , denoted $\boldsymbol{<_{\tau,i}}$, is the monomial order in $R_{n,i}$ with the property:

{\centering $\underline{\alpha}=(\alpha_1, \dots, \alpha_{i-1},0,\alpha_{i+1}, \dots, \alpha_{n}) <_{\tau} (\beta_1, \dots, \beta_{i-1},0,\beta_{i+1}, \dots, \beta_{n})=\underline{\beta}$ 

if and only if 

$(\alpha_1, \dots, \alpha_{i-1},\alpha_{i+1}, \dots, \alpha_{n}) <_{\tau,i} (\beta_1, \dots, \beta_{i-1},\beta_{i+1}, \dots, \beta_{n})$.\par}}

From this point on instead of $(\alpha_1, \dots, \alpha_{i-1},\alpha_{i+1}, \dots, \alpha_{n}) <_{\tau,i} (\beta_1, \dots, \beta_{i-1},\beta_{i+1}, \dots, \beta_{n})$ we will write $(\alpha_1, \dots, \alpha_{i-1},0,\alpha_{i+1}, \dots, \alpha_{n}) <_{\tau,i} (\beta_1, \dots, \beta_{i-1},0,\beta_{i+1}, \dots, \beta_{n})$.

\

The definition above allows us to state the following theorem.

{\theorem \label{lexprop} If $n>3$ and $<_{\tau,i}$ is the lexicographic (resp. graded lexicographic, reverse lexicographic) order in $R_{n,i}$ for all $1 \leq i \leq n$ then $<_{\tau}$ is the lexicographic (resp. graded lexicographic, reverse lexicographic) in $R_n$. \par}

\begin{proof}

We will first prove the case for the lexicographic order. Assume that $\underline{\alpha} <_{lex} \underline{\beta}$, and let $k=\text{min}\{j: \alpha_j \neq \beta_j \}$.

\begin{itemize}

\item If $k>1$ then $(\alpha_1, \dots, \alpha_{k-1}, \alpha_k, \dots, \alpha_n)<_{\tau} (\beta_1, \dots, \beta_{k-1}, \beta_k, \dots, \beta_n)$ if and only if $(\alpha_1, \dots, \alpha_{k-1}, \alpha_k, \dots, \alpha_n)<_{lex} (\beta_1, \dots, \beta_{k-1}, \beta_k, \dots, \beta_n)$ since $(0,\dots,\alpha_k, \dots, \alpha_n) <_{lex,1} (0, \dots,\beta_k, \dots, \beta_n)$ and $<_{\tau,1}=<_{lex,1}$.

\item If $k=1$ then $(\alpha_1, \dots, \alpha_n) <_{\tau} (\beta_1, \dots, \beta_n)$ if and only if $(\alpha_1, \dots, \alpha_n) <_{lex} (\beta_1, \dots, \beta_n)$, due to the sequence $(\alpha_1,\alpha_2, \dots, \alpha_n) \leq_{lex,1} (\alpha_1, \sum_{i=2}^{n} \alpha_i, \dots, 0) <_{lex,n} (\beta_1, \dots, 0) \leq_{lex,1} (\beta_1, \dots, \beta_n)$ and the fact that $<_{\tau,i}=<_{lex,i}$ for all $1 \leq i \leq n$.

\end{itemize}

For the case of the reverse lexicographic order notice that if there is $i$ such that $\alpha_{i}=\beta_{i}$ then $\underline{\alpha} <_{\tau} \underline{\beta}$  if and only if  $(\alpha_1, \dots, \alpha_{i-1},0,\alpha_{i+1}, \dots, \alpha_n)<_{\tau} (\beta_1, \dots, \beta_{i-1},0,\beta_{i+1}, \dots, \beta_n)$, which is equivalent to saying $(\alpha_1, \dots, \alpha_{i-1}, \alpha_{i+1},\dots,\alpha_n) <_{\tau,i} (\beta_1,\dots, \beta_{i-1}, \beta_{i+1},\dots, \beta_n)$, and this is equivalent to $(\alpha_1, \dots, \alpha_{i-1},0,\alpha_{i+1}, \dots, \alpha_n)<_{revlex} (\beta_1, \dots, \beta_{i-1},0,\beta_{i+1}, \dots, \beta_n)$, since $<_{\tau,i}=<_{revlex,i}$, and this happens if and only if $\underline{\alpha} <_{revlex} \underline{\beta}$. 

Hence $(d,0,\dots,0)<_{\tau} (0,0,\dots, d+1)$ if and only if $(d,0,\dots,0) <_{revlex} (0,0,\dots,d+1)$.

Additionally if $\gamma_i=\min{(\alpha_i,\beta_i)}$ then $\underline{\alpha}<\underline{\beta}$ if and only if $\underline{\alpha}-\underline{\gamma}<\underline{\beta}-\underline{\gamma}$ for any monomial order $<$.

Because of the observations above, to guarantee that $<_{\tau}=<_{revlex}$ it is enough to prove that $\underline{\alpha}<_{\tau} \underline{\beta}$ if and only if $\underline{\alpha} <_{revlex} \underline{\beta}$ for all $\underline{\alpha},\underline{\beta}$ such that $\sum_{i=1}^n \alpha_i=\sum_{i=1}^n \beta_i=d$, $\alpha_i \cdot \beta_i=0$ and $\alpha_i+\beta_i>0$. So we will restrict ourselves to these cases from this point on. 

Let's assume without loss of generality that $\alpha_n \neq 0$, and let $k=\max \{j: \beta_j \neq 0 \}$, so $k<n$. We will divide the argument in five cases:

\begin{itemize}
\item \textit{Case 1}: $k<n-1$

Notice that $(\alpha_1, \dots,0, \alpha_{k+1}, \dots, \alpha_{n-1}, \alpha_n)<_{revlex} (\beta_1,\dots, \beta_k, 0, \dots, 0, 0)$ if and only if $(\alpha_1, \dots,0, \alpha_{k+1}, \dots, \alpha_{n-1}, \alpha_n)<_{\tau} (\beta_1,\dots, \beta_k, 0, \dots, 0, 0)$, due to the sequence:

{\centering $(\alpha_1, \dots,0, \alpha_{k+1}, \dots, \alpha_{n-1}, \alpha_n)<_{revlex,k} (0, \dots,0,0, \dots,d,0)$ 

$(0, \dots,0,0, \dots,d,0)<_{revlex,n}  (0, \dots,d,0,\dots, 0, 0) \leq_{revlex,n} (\beta_1,\dots, \beta_k, 0, \dots, 0, 0)$ \par}

and the fact that $<_{\tau,i}=<_{revlex,i}$ for all $1 \leq i \leq n$.

\item \textit{Case 2}: $k=n-1$ and $\alpha_i \neq 0$ for all $i<n-1$

Notice that $(\alpha_1, \dots, \alpha_{n-2},0,\alpha_n)<_{revlex} (0, \dots, 0,d,0)$ if and only if $(\alpha_1, \dots, \alpha_{n-2},0,\alpha_n)<_{\tau} (0, \dots, 0,d,0)$, due to the sequence:

{\centering $(\alpha_1, \dots, \alpha_{n-2}, 0, \alpha_n)<_{revlex,n-1} (\alpha_1+\alpha_{n-2}, \dots,0,0,\alpha_n)<_{revlex,n-2}(0,\dots,0,d,0)$ 
 \par}

and the fact that $<_{\tau,i}=<_{revlex,i}$ for all $1 \leq i \leq n$.

\item \textit{Case 3}: $k=n-1$, $\alpha_{i}=0$ for all $i<n-1$.

Notice that $(0,\dots,0,0,d)<_{revlex} (\beta_1,\dots,\beta_{n-2},\beta_{n-1},0)$ if and only if $(0,\dots,0,0,d)<_{\tau} (\beta_1,\dots,\beta_{n-2},\beta_{n-1},0)$, due to the sequence:
 
{\centering $(0,\dots,0,0,d)<_{revlex,1} (0,\dots,0,d,0) <_{revlex,n} (\beta_1,\dots,\beta_{n-2},\beta_{n-1},0)$ 
 \par}

and the fact that $<_{\tau,i}=<_{revlex,i}$ for all $1 \leq i \leq n$.

\item \textit{Case 4}: $k=n-1$, $\alpha_{n-2}=0$ and there is $l<n-2$ with $\alpha_l\neq0$.

Notice that $(\alpha_1,\dots,\alpha_l, 0, \dots, 0,0,\alpha_n)<_{revlex} (\beta_1,\dots,0,\beta_{l+1},\dots,\beta_{n-2},\beta_{n-1},0)$ if and only if $(\alpha_1,\dots,\alpha_l, 0, \dots, 0,0,\alpha_n)<_{\tau} (\beta_1,\dots,0,\beta_{l+1},\dots,\beta_{n-2},\beta_{n-1},0)$, due to the sequence:
 
{\centering $(\alpha_1,\dots,\alpha_l, 0, \dots, 0,0,\alpha_n)<_{revlex,n-2} (0,\dots,0,d,0) <_{revlex,n} (\beta_1,\dots,0,\beta_{l+1},\dots,\beta_{n-2},\beta_{n-1},0)$ 
 \par}

and the fact that $<_{\tau,i}=<_{revlex,i}$ for all $1 \leq i \leq n$.

\item \textit{Case 5:} $k=n-1$, $\alpha_{n-2}\neq 0$ and there is $l<n-2$ with $\alpha_l =0$.

Notice that $(\alpha_1,\dots,0,\dots,\alpha_{n-2},0,\alpha_n)<_{revlex} (\beta_1,\dots,\beta_{l},\dots,0,\beta_{n-1},0)$ if and only if $(\alpha_1,\dots,0,\dots,\alpha_{n-2},0,\alpha_n)<_{\tau} (\beta_1,\dots,\beta_{l},\dots,0,\beta_{n-1},0)$, due to the sequence:
 
{\centering $(\alpha_1,\dots,0,\dots,\alpha_{n-2},0,\alpha_n)<_{revlex,l} (0,\dots,0,\dots,0,d,0) <_{revlex,n} (\beta_1,\dots,\beta_{l},\dots,0,\beta_{n-1},0)$ 
 \par}

and the fact that $<_{\tau,i}=<_{revlex,i}$ for all $1 \leq i \leq n$.
\end{itemize}

The case for the graded lexicographic order is analogous to that of the reverse lexicographic order.

This completes the proof.
\end{proof}

Theorem \ref{lexprop},  coupled with Remark \ref{spread}, motivated A. Conca, \cite{Conca}, to ask the following:

{\question \label{question} Is there $n > 3$ such that any monomial ordering in $R_n$ is uniquely determined by its induced orderings? \par}

\

The answer to this question is negative, and its proof is given in Theorem \ref{main}. Furthermore the answer is still negative if we just focus our attention on graded monomial orders.

{\theorem \label{matrix} \textbf{[Robbiano]} Given a monomial order $<_{\tau}$ in $R_n$ there exists $A \in M_{m \times n}(\mathbb{R})$ such that $\underline{\alpha_1} <_{\tau} \underline{\alpha_2}$ if and only if $A\cdot \underline{\alpha_1} <_{lex} A \cdot \underline{\alpha_2}$. \par}

Proofs of this result can be found in \cite{GR}, \cite{OZ}, and \cite{Ro}.

{\example \label{matrices}  

The $n \times n$ identity matrix $I$ corresponds to the the $<_{lex}$ order in $R_n$. While the matrices $G = \begin{bmatrix} 1 & 1 & \cdots & 1 & 1 & 1 \\ 1 & 0 & \cdots & 0 & 0 &0 \\ 0 & 1 & \cdots & 0 & 0 & 0 \\ \vdots & \vdots & \ddots& \vdots & \vdots & \vdots \\ 0 & 0 & \cdots & 1 & 0 & 0 \\ 0 & 0 & \cdots & 0 & 1 & 0 \end{bmatrix}$ and $R=\begin{bmatrix*}[r] 1 & 1 & 1 & \cdots & 1 & 1  \\ 0 & 0 & 0 & \cdots  & 0 & -1 \\  0 & 0 & 0 & \cdots & -1 &0 \\ \vdots & \vdots & \vdots & \ddots & \vdots & \vdots  \\  0 & 0 & -1 & \cdots & 0 & 0 \\ 0 & -1 & 0 & \cdots & 0 & 0  \end{bmatrix*} $ correspond to $<_{deglex}$ and  $<_{revlex}$ respectively.
\par}

We provide now a converse of Theorem \ref{matrix}, which is presented as a statement that encompasses \cite[Exercise 2.8]{CLO2} and \cite[Proposition 1.4.12]{KR}

{\lemma \label{nonzero} Let $A \in M_{m \times n}(\mathbb{R})$ such that $\text{Ker}(A) \cap \mathbb{Z}^n=\{(0,\dots,0)\}$ and the first non-zero entry in each column is positive. Then there is a monomial order $<_{A}$ in $R_n$ with $\underline{\alpha_1} <_{A} \underline{\alpha_2}$ if $A\cdot \underline{\alpha_1} <_{lex} A \cdot \underline{\alpha_2}$. \par}

\

In particular if $A \in M_{n \times n}(\mathbb{N}^{*})$ and $\text{det}(A) \neq 0$ then $<_{A}$ is a monomial order.

{\lemma \label{lower} Let $A,B  \in M_{m \times n}(\mathbb{R})$ be matrices defining monomial orders and let $L \in M_{m \times m}(\mathbb{R})$ be a lower triangular matrix with positive entries in the diagonal such that $B=LA$ . Then the monomial orders defined by $A$ and $B$ are the same. \par}

\

An elementary proof of Lemma \ref{lower} appears in \cite{OZ}, it is also an exercise \cite[Tutorial 9]{KR}.

\

Additionally if we have a matrix representation for a monomial order $<_{\tau}$ in terms of an $n \times n$ matrix with entries on the integers, we can find the matrix representation for its $i^{th}$  induced ordering because of the following proposition in \cite[Proposition 1.4.13]{KR}.

{\lemma \label{induced} If $A$ is an $n \times n$ matrix representing the monomial order $<_{\tau}$ then its $i^{th}$ induced ordering $<_{\tau,i}$ is represented by the matrix $A_i$ which is obtained by first deleting the $i^{th}$ column of A and then the first row which is linearly independent on those above it. \par}

\

We are now ready to prove our main result.

{\theorem \label{main} For $n \geq 4$ there exist $<_{\tau}$ and $<_{\tau'}$ distinct monomial orders in $R_n$ such that their induced orderings $<_{\tau,i}$ and $<_{\tau',i}$ are the same in $R_{n,i}$ for all $1 \leq i \leq n$. \par}

\begin{proof}

Consider the $n \times n$ matrices, $n \geq 4$, 

{\centering $C_n=\begin{bmatrix} 1 & 1 & 1 & \cdots & 1 & 1 & 1 & 1 \\ 1 & 1& 0 & \cdots & 0 & 0 & 0 &0 \\ 1 & 0 & 1 & \cdots & 0 & 0 & 0 & 0 \\ \vdots & \vdots & \vdots & \ddots & \vdots & \vdots & \vdots & \vdots \\ 1 & 0 & 0 & \cdots & 0 & 1 & 0 & 0  \\ \displaystyle \frac{n^2+n+2}{2} & n-1 & n-2 & \cdots & 4 & 3 & 2 & 1 \\ 0 & 0 & 0 & \cdots & 0 & 0 & 0 & 1 \end{bmatrix}$ \par}

and

{\centering $D_n=\begin{bmatrix} 1 & 1 & 1 & \cdots & 1 & 1 & 1 & 1 \\ 1 & 1& 0 & \cdots & 0 & 0 & 0 &0 \\ 1 & 0 & 1 & \cdots & 0 & 0 & 0 & 0 \\ \vdots & \vdots & \vdots & \ddots & \vdots & \vdots & \vdots & \vdots \\ 1 & 0 & 0 & \cdots & 0 & 1 & 0 & 0  \\ \displaystyle \frac{n^2-n}{2} & n-1 & n-2 & \cdots & 4 & 2 & 3 & 1 \\ 0 & 0 & 0 & \cdots & 0 & 0 & 0 & 1 \end{bmatrix}$. \par}

A simple reduction allows us to calculate $\text{det}(C_n)=4-3n$ and $\text{det}(D_n)=5-2n$. Since $4-3n\neq 0 \neq 5-2n$ for any value of $n>3$, Lemma \ref{nonzero} proves that they both define monomial orders on $R_n$.

Notice that the monomial orders defined by $C_n$ and $D_n$ are distinct in $R_n$, since  $(2,0, \dots, 0, n, n^2,2) <_{C_n} (4,0, \dots, 0, n^2,n,1)$,$(4,0, \dots, 0, n^2, n,1) <_{D_n} (2,0, \dots, 0,n, n^2, 2)$.

Let $C_{n,i}$ (resp. $D_{n,i}$) be the $(n-1) \times (n-1)$ matrix obtained by eliminating the $i^{th}$ column and the $n^{th}$ row of $C_n$ (resp. $D_n$). We will now present the numerical values of the determinants of $C_{n,i}$ and $D_{n,i}$:  

\begin{multicols}{2}
$\text{det}(C_{n,1}) = (-1)^n$

$\text{det}(C_{n,i}) =  (-1)^{n+i-2} $

$\text{det}(C_{n,n-1}) = -2n$ 

$\text{det}(C_{n,n}) = 4-3n$ 

$\text{det}(D_{n,1}) = (-1)^n\cdot 2$ 

$\text{det}(D_{n,i}) =  (-1)^{n+i-2}\cdot 2$   

$\text{det}(D_{n,n-1}) =  -1-2n$  

$\text{det}(D_{n,n}) =  5-2n$ 
\end{multicols}
\noindent for $ 2 \leq i \leq n-2$.

All these determinants are nonzero, for $n>3$, so Lemma \ref{induced} implies that $C_{n,i}$ (resp. $D_{n,i}$) are matrix representations for the $\text{i}^{th}$ induced ordering of $<_{C_n}$ (resp. $<_{D_n}$) in $R_{n,i}$.

Finally the facts that $C_{n,i}$ and $D_{n,i}$ are invertible and their first $n-2$ rows are the same guarantee that there is a lower triangular matrix 

{\centering $U_i=\begin{bmatrix} 1 & 0 & \cdots & 0 & 0 \\ 0 & 1 & \cdots & 0 & 0 \\ \vdots & \vdots& \ddots & \vdots & \vdots \\ 0 & 0 & \cdots & 1 & 0 \\ a_{i,1} & a_{i,2} & \cdots & a_{i,n-2} & a_{i,n-1}  \end{bmatrix}$, \par} 

\noindent such that $D_{n,i}=U_iC_{n,i}$. By Cramer's rule $\displaystyle a_{i,n-1}=\frac{\text{det}(D_{n,i})}{\text{det}(C_{n,i})}>0$ which implies by Lemma \ref{lower}, that the induced orders $C_{n,i}$ and $D_{n,i}$ are the same for all $1 \leq i \leq n$. 

This concludes our proof.
\end{proof}
\subsection*{Acknowledgements}

I would like to thank the organizers of the USTARS conference for inviting me to present at both USTARS 2013, where I gave a talk on the main result of this paper, and USTARS 2014. I would also like to thank Prof. Alejandra Alvarado and Prof. Edray Goins for encouraging me to attend and present at USTARS 2013 and 2014, and to Prof. Giulio Caviglia for bringing to my attention Question \ref{question}.


\bibliographystyle{plain}
\bibliography{refPRF}

\end{document}